\documentclass{gtart}

\def\ifplaintex{\expandafter\ifx\csname documentclass\endcsname\relax}

\def\gtp{{\mathsurround=0pt\it $\cal G\mskip-2mu$eometry \&\ 
$\cal T\!\!$opology $\cal P\!$ublications}}  

\def\recd{{\small Received:\qua\receiveddate\ifx\reviseddate\relax
\else\qquad Revised:\qua\reviseddate\fi\par}} 

\def\lognumber#1{\def\thelognumber{#1}}
\def\volumenumber#1{\def\thevolumenumber{#1}}
\def\volumeyear#1{\def\thevolumeyear{#1}}
\def\papernumber#1{\def\thepapernumber{#1}}
\def\pagenumbers#1#2{\def\startpage{#1}\def\finishpage{#2}}
\def\published#1{\def\publishdate{#1}}

\def\received#1{\def\receiveddate{#1}}

\def\accepted#1{\def\accepteddate{#1}}

\def\covertitle#1{\def\thecovertitle{#1}}

\def\asciiaddress#1{\def\theasciiaddress{#1}}

\long\def\asciiabstract#1{\long\def\theasciiabstract{#1}}

\let\\\par\let\thelognumber\relax\let\thevolumenumber\relax
\let\thepapernumber\relax\let\thevolumeyear\relax\let\startpage\relax
\let\finishpage\relax\let\publishdate\relax\let\receiveddate\relax
\let\reviseddate\relax\let\accepteddate\relax\let\theasciititle\relax
\let\thecovertitle\relax\let\theasciiauthors\relax\let\theasciiaddress\relax
\let\theasciiabstract\relax

\let\theasciiemail\relax

\ifplaintex
\font\logobig=cmssbx10 scaled 3836
\font\logomed=cmssbx10 scaled 2557
\else
\font\logobig=cmssbx10 scaled 4200
\font\logomed=cmssbx10 scaled 2800
\fi

\long\def\makeagttitle{   
\count0=\startpage
\agt\hfill      
\hbox to 45truept{\vbox to 0pt{\vglue -13truept{\logomed A\kern -.37em{\logobig 
T}\kern -.38em G}\vss}\hss}
\break
{\small Volume \thevolumenumber\ (\thevolumeyear)
\startpage--\finishpage\nl
Published: \publishdate}

\vglue .25truein

{\parskip=0pt\leftskip 0pt plus
1fil\def\\{\par\smallskip}{\Large\bf\thetitle}\par\medskip} \vglue
0.05truein

{\parskip=0pt\leftskip 0pt plus 1fil\def\\{\par}{\sc\theauthors}
\par\medskip}%
 
\vglue 0.03truein 

{\small\leftskip 25truept\rightskip 25truept{\bf Abstract}\stdspace\theabstract

{\bf AMS Classification}\stdspace\theprimaryclass
\ifx\thesecondaryclass\relax\else; \thesecondaryclass\fi\par
{\bf Keywords}\stdspace \thekeywords\par}\vglue 7truept

}   

\ifplaintex
\font\phead=cmsl9 scaled 950
\font\pnum=cmbx10 scaled 913
\font\pfoot=cmsl9 scaled 950
\headline{\vbox to 0pt{\vskip -4.5mm\line{\small\phead\ifnum
\count0=\startpage ISSN 1472-2739 (on-line) 1472-2747 (printed)
\hfill {\pnum\folio}\else\ifodd\count0\def\\{ }%
\ifx\theshorttitle\relax\thetitle\else\theshorttitle\fi\hfill{\pnum\folio}
\else\def\\{ and }{\pnum\folio}\hfill\ifx\theshortauthors\relax\theauthors
\else\theshortauthors\fi\fi\fi}\vss}}
\footline{\vbox to 0pt{\vglue 0mm\line{\small\pfoot\ifnum\count0=\startpage
\copyright\ \gtp\hfill\else
\agt, Volume \thevolumenumber\ (\thevolumeyear)\hfill\fi}\vss}}
\else
\font=cmsl9 scaled 1050
\font\lnum=cmbx10 
\font=cmsl9 scaled 1050
\makeatletter
\def\@oddhead{{\small\lhead\ifnum\count0=\startpage ISSN 1472-2739 
(on-line) 1472-2747 (printed)\hfill {\lnum\number\count0}\else\ifodd\count0
\def\\{ }\ifx\theshorttitle\relax \thetitle \else\theshorttitle\fi\hfill
{\lnum\number\count0}\else\def\\{ and }{\lnum\number\count0}
\hfill\ifx\theshortauthors\relax 
\theauthors\else\theshortauthors\fi\fi\fi}}\def\@evenhead{\@oddhead}
\def\@oddfoot{\small\lfoot\ifnum\count0=\startpage\copyright\ \gtp\hfill\else
\agt, Volume \thevolumenumber\ (\thevolumeyear)\hfill\fi}
\def\@evenfoot{\@oddfoot}
\makeatother
\fi
\let\maketitlepage\makeagttitle
\let\makeshorttitle\maketitlepage
\let\maketitle\maketitlepage

\newwrite\gtoutfile
\long\gdef\makeheadfile{  
{\def\\{, }\def\s{ }
\immediate\openout\gtoutfile head.xxx
\immediate\write\gtoutfile{To: math@arxiv.org}
\immediate\write\gtoutfile{Subject: put OR rep NNNNN:ppppp}
\immediate\write\gtoutfile{--text follows this line--}
\immediate\write\gtoutfile{Proxy-for: \ifx\theasciiauthors\relax
\theauthors\else\theasciiauthors\fi\s<\ifx\theasciiemail\relax\theemail\else\theasciiemail\fi>}
\immediate\write\gtoutfile{\noexpand\\}
\immediate\write\gtoutfile{Authors: \ifx\theasciiauthors\relax
\theauthors\else\theasciiauthors\fi}
{\def\\{ }\immediate\write\gtoutfile{Title: \ifx\theasciititle\relax
\thetitle\else\theasciititle\fi}}
\immediate\write\gtoutfile{Subj-class: GT or SG, GR etc}
\immediate\write\gtoutfile{MSC-class: \theprimaryclass\ifx\thesecondaryclass\relax\else, \thesecondaryclass\fi}
\immediate\write\gtoutfile{Journal-ref: Algebr. Geom. Topol. \thevolumenumber\s
(\thevolumeyear) \startpage-\finishpage}
\immediate\write\gtoutfile{Comments: Published by Algebraic and
Geometric Topology at}
\immediate\write\gtoutfile{\s\s\s  http://www.maths.warwick.ac.uk/agt/AGTVol\thevolumenumber/agt-\thevolumenumber-\thepapernumber.abs.html}
\immediate\write\gtoutfile{\noexpand\\}
\immediate\write\gtoutfile{}
\ifx\theasciiabstract\relax
\immediate\write\gtoutfile{\theabstract}\else
\immediate\write\gtoutfile{\theasciiabstract}\fi
\immediate\write\gtoutfile{}
\immediate\write\gtoutfile{\noexpand\\}
\immediate\write\gtoutfile{}
\immediate\closeout\gtoutfile}}  

\def\maketitlepage{\makeagttitle\makeheadfile}
\let\makeshorttitle\maketitlepage
\let\maketitle\maketitlepage

\def\ifplaintex{\expandafter\ifx\csname documentclass\endcsname\relax}

\def\gtp{{\mathsurround=0pt\it $\cal G\mskip-2mu$eometry \&\ 
$\cal T\!\!$opology $\cal P\!$ublications}}  

\def\recd{{\small Received:\qua\receiveddate\ifx\reviseddate\relax
\else\qquad Revised:\qua\reviseddate\fi\par}} 

\def\lognumber#1{\def\thelognumber{#1}}
\def\volumenumber#1{\def\thevolumenumber{#1}}
\def\volumeyear#1{\def\thevolumeyear{#1}}
\def\papernumber#1{\def\thepapernumber{#1}}
\def\pagenumbers#1#2{\def\startpage{#1}\def\finishpage{#2}}
\def\published#1{\def\publishdate{#1}}

\def\received#1{\def\receiveddate{#1}}

\def\accepted#1{\def\accepteddate{#1}}

\def\covertitle#1{\def\thecovertitle{#1}}

\def\asciiaddress#1{\def\theasciiaddress{#1}}

\long\def\asciiabstract#1{\long\def\theasciiabstract{#1}}

\let\\\par\let\thelognumber\relax\let\thevolumenumber\relax
\let\thepapernumber\relax\let\thevolumeyear\relax\let\startpage\relax
\let\finishpage\relax\let\publishdate\relax\let\receiveddate\relax
\let\reviseddate\relax\let\accepteddate\relax\let\theasciititle\relax
\let\thecovertitle\relax\let\theasciiauthors\relax\let\theasciiaddress\relax
\let\theasciiabstract\relax

\let\theasciiemail\relax

\ifplaintex
\font\logobig=cmssbx10 scaled 3836
\font\logomed=cmssbx10 scaled 2557
\else
\font\logobig=cmssbx10 scaled 4200
\font\logomed=cmssbx10 scaled 2800
\fi

\long\def\makeagttitle{   
\count0=\startpage
\agt\hfill      
\hbox to 45truept{\vbox to 0pt{\vglue -13truept{\logomed A\kern -.37em{\logobig 
T}\kern -.38em G}\vss}\hss}
\break
{\small Volume \thevolumenumber\ (\thevolumeyear)
\startpage--\finishpage\nl
Published: \publishdate}

\vglue .25truein

{\parskip=0pt\leftskip 0pt plus
1fil\def\\{\par\smallskip}{\Large\bf\thetitle}\par\medskip} \vglue
0.05truein

{\parskip=0pt\leftskip 0pt plus 1fil\def\\{\par}{\sc\theauthors}
\par\medskip}%
 
\vglue 0.03truein 

{\small\leftskip 25truept\rightskip 25truept{\bf Abstract}\stdspace\theabstract

{\bf AMS Classification}\stdspace\theprimaryclass
\ifx\thesecondaryclass\relax\else; \thesecondaryclass\fi\par
{\bf Keywords}\stdspace \thekeywords\par}\vglue 7truept

}   

\ifplaintex
\font\phead=cmsl9 scaled 950
\font\pnum=cmbx10 scaled 913
\font\pfoot=cmsl9 scaled 950
\headline{\vbox to 0pt{\vskip -4.5mm\line{\small\phead\ifnum
\count0=\startpage ISSN 1472-2739 (on-line) 1472-2747 (printed)
\hfill {\pnum\folio}\else\ifodd\count0\def\\{ }%
\ifx\theshorttitle\relax\thetitle\else\theshorttitle\fi\hfill{\pnum\folio}
\else\def\\{ and }{\pnum\folio}\hfill\ifx\theshortauthors\relax\theauthors
\else\theshortauthors\fi\fi\fi}\vss}}
\footline{\vbox to 0pt{\vglue 0mm\line{\small\pfoot\ifnum\count0=\startpage
\copyright\ \gtp\hfill\else
\agt, Volume \thevolumenumber\ (\thevolumeyear)\hfill\fi}\vss}}
\else
\font=cmsl9 scaled 1050
\font\lnum=cmbx10 
\font=cmsl9 scaled 1050
\makeatletter
\def\@oddhead{{\small\lhead\ifnum\count0=\startpage ISSN 1472-2739 
(on-line) 1472-2747 (printed)\hfill {\lnum\number\count0}\else\ifodd\count0
\def\\{ }\ifx\theshorttitle\relax \thetitle \else\theshorttitle\fi\hfill
{\lnum\number\count0}\else\def\\{ and }{\lnum\number\count0}
\hfill\ifx\theshortauthors\relax 
\theauthors\else\theshortauthors\fi\fi\fi}}\def\@evenhead{\@oddhead}
\def\@oddfoot{\small\lfoot\ifnum\count0=\startpage\copyright\ \gtp\hfill\else
\agt, Volume \thevolumenumber\ (\thevolumeyear)\hfill\fi}
\def\@evenfoot{\@oddfoot}
\makeatother
\fi
\let\maketitlepage\makeagttitle
\let\makeshorttitle\maketitlepage
\let\maketitle\maketitlepage

\newwrite\gtoutfile
\long\gdef\makeheadfile{  
{\def\\{, }\def\s{ }
\immediate\openout\gtoutfile head.xxx
\immediate\write\gtoutfile{To: math@arxiv.org}
\immediate\write\gtoutfile{Subject: put OR rep NNNNN:ppppp}
\immediate\write\gtoutfile{--text follows this line--}
\immediate\write\gtoutfile{Proxy-for: \ifx\theasciiauthors\relax
\theauthors\else\theasciiauthors\fi\s<\ifx\theasciiemail\relax\theemail\else\theasciiemail\fi>}
\immediate\write\gtoutfile{\noexpand\\}
\immediate\write\gtoutfile{Authors: \ifx\theasciiauthors\relax
\theauthors\else\theasciiauthors\fi}
{\def\\{ }\immediate\write\gtoutfile{Title: \ifx\theasciititle\relax
\thetitle\else\theasciititle\fi}}
\immediate\write\gtoutfile{Subj-class: GT or SG, GR etc}
\immediate\write\gtoutfile{MSC-class: \theprimaryclass\ifx\thesecondaryclass\relax\else, \thesecondaryclass\fi}
\immediate\write\gtoutfile{Journal-ref: Algebr. Geom. Topol. \thevolumenumber\s
(\thevolumeyear) \startpage-\finishpage}
\immediate\write\gtoutfile{Comments: Published by Algebraic and
Geometric Topology at}
\immediate\write\gtoutfile{\s\s\s  http://www.maths.warwick.ac.uk/agt/AGTVol\thevolumenumber/agt-\thevolumenumber-\thepapernumber.abs.html}
\immediate\write\gtoutfile{\noexpand\\}
\immediate\write\gtoutfile{}
\ifx\theasciiabstract\relax
\immediate\write\gtoutfile{\theabstract}\else
\immediate\write\gtoutfile{\theasciiabstract}\fi
\immediate\write\gtoutfile{}
\immediate\write\gtoutfile{\noexpand\\}
\immediate\write\gtoutfile{}
\immediate\closeout\gtoutfile}}  

\def\maketitlepage{\makeagttitle\makeheadfile}
\let\makeshorttitle\maketitlepage
\let\maketitle\maketitlepage

\def\ifplaintex{\expandafter\ifx\csname documentclass\endcsname\relax}

\def\gtp{{\mathsurround=0pt\it $\cal G\mskip-2mu$eometry \&\ 
$\cal T\!\!$opology $\cal P\!$ublications}}  

\def\recd{{\small Received:\qua\receiveddate\ifx\reviseddate\relax
\else\qquad Revised:\qua\reviseddate\fi\par}} 

\def\lognumber#1{\def\thelognumber{#1}}
\def\volumenumber#1{\def\thevolumenumber{#1}}
\def\volumeyear#1{\def\thevolumeyear{#1}}
\def\papernumber#1{\def\thepapernumber{#1}}
\def\pagenumbers#1#2{\def\startpage{#1}\def\finishpage{#2}}
\def\published#1{\def\publishdate{#1}}

\def\received#1{\def\receiveddate{#1}}

\def\accepted#1{\def\accepteddate{#1}}

\def\covertitle#1{\def\thecovertitle{#1}}

\def\asciiaddress#1{\def\theasciiaddress{#1}}

\long\def\asciiabstract#1{\long\def\theasciiabstract{#1}}

\let\\\par\let\thelognumber\relax\let\thevolumenumber\relax
\let\thepapernumber\relax\let\thevolumeyear\relax\let\startpage\relax
\let\finishpage\relax\let\publishdate\relax\let\receiveddate\relax
\let\reviseddate\relax\let\accepteddate\relax\let\theasciititle\relax
\let\thecovertitle\relax\let\theasciiauthors\relax\let\theasciiaddress\relax
\let\theasciiabstract\relax

\let\theasciiemail\relax

\ifplaintex
\font\logobig=cmssbx10 scaled 3836
\font\logomed=cmssbx10 scaled 2557
\else
\font\logobig=cmssbx10 scaled 4200
\font\logomed=cmssbx10 scaled 2800
\fi

\long\def\makeagttitle{   
\count0=\startpage
\agt\hfill      
\hbox to 45truept{\vbox to 0pt{\vglue -13truept{\logomed A\kern -.37em{\logobig 
T}\kern -.38em G}\vss}\hss}
\break
{\small Volume \thevolumenumber\ (\thevolumeyear)
\startpage--\finishpage\nl
Published: \publishdate}

\vglue .25truein

{\parskip=0pt\leftskip 0pt plus
1fil\def\\{\par\smallskip}{\Large\bf\thetitle}\par\medskip} \vglue
0.05truein

{\parskip=0pt\leftskip 0pt plus 1fil\def\\{\par}{\sc\theauthors}
\par\medskip}%
 
\vglue 0.03truein 

{\small\leftskip 25truept\rightskip 25truept{\bf Abstract}\stdspace\theabstract

{\bf AMS Classification}\stdspace\theprimaryclass
\ifx\thesecondaryclass\relax\else; \thesecondaryclass\fi\par
{\bf Keywords}\stdspace \thekeywords\par}\vglue 7truept

}   

\ifplaintex
\font\phead=cmsl9 scaled 950
\font\pnum=cmbx10 scaled 913
\font\pfoot=cmsl9 scaled 950
\headline{\vbox to 0pt{\vskip -4.5mm\line{\small\phead\ifnum
\count0=\startpage ISSN 1472-2739 (on-line) 1472-2747 (printed)
\hfill {\pnum\folio}\else\ifodd\count0\def\\{ }%
\ifx\theshorttitle\relax\thetitle\else\theshorttitle\fi\hfill{\pnum\folio}
\else\def\\{ and }{\pnum\folio}\hfill\ifx\theshortauthors\relax\theauthors
\else\theshortauthors\fi\fi\fi}\vss}}
\footline{\vbox to 0pt{\vglue 0mm\line{\small\pfoot\ifnum\count0=\startpage
\copyright\ \gtp\hfill\else
\agt, Volume \thevolumenumber\ (\thevolumeyear)\hfill\fi}\vss}}
\else
\font=cmsl9 scaled 1050
\font\lnum=cmbx10 
\font=cmsl9 scaled 1050
\makeatletter
\def\@oddhead{{\small\lhead\ifnum\count0=\startpage ISSN 1472-2739 
(on-line) 1472-2747 (printed)\hfill {\lnum\number\count0}\else\ifodd\count0
\def\\{ }\ifx\theshorttitle\relax \thetitle \else\theshorttitle\fi\hfill
{\lnum\number\count0}\else\def\\{ and }{\lnum\number\count0}
\hfill\ifx\theshortauthors\relax 
\theauthors\else\theshortauthors\fi\fi\fi}}\def\@evenhead{\@oddhead}
\def\@oddfoot{\small\lfoot\ifnum\count0=\startpage\copyright\ \gtp\hfill\else
\agt, Volume \thevolumenumber\ (\thevolumeyear)\hfill\fi}
\def\@evenfoot{\@oddfoot}
\makeatother
\fi
\let\maketitlepage\makeagttitle
\let\makeshorttitle\maketitlepage
\let\maketitle\maketitlepage

\newwrite\gtoutfile
\long\gdef\makeheadfile{  
{\def\\{, }\def\s{ }
\immediate\openout\gtoutfile head.xxx
\immediate\write\gtoutfile{To: math@arxiv.org}
\immediate\write\gtoutfile{Subject: put OR rep NNNNN:ppppp}
\immediate\write\gtoutfile{--text follows this line--}
\immediate\write\gtoutfile{Proxy-for: \ifx\theasciiauthors\relax
\theauthors\else\theasciiauthors\fi\s<\ifx\theasciiemail\relax\theemail\else\theasciiemail\fi>}
\immediate\write\gtoutfile{\noexpand\\}
\immediate\write\gtoutfile{Authors: \ifx\theasciiauthors\relax
\theauthors\else\theasciiauthors\fi}
{\def\\{ }\immediate\write\gtoutfile{Title: \ifx\theasciititle\relax
\thetitle\else\theasciititle\fi}}
\immediate\write\gtoutfile{Subj-class: GT or SG, GR etc}
\immediate\write\gtoutfile{MSC-class: \theprimaryclass\ifx\thesecondaryclass\relax\else, \thesecondaryclass\fi}
\immediate\write\gtoutfile{Journal-ref: Algebr. Geom. Topol. \thevolumenumber\s
(\thevolumeyear) \startpage-\finishpage}
\immediate\write\gtoutfile{Comments: Published by Algebraic and
Geometric Topology at}
\immediate\write\gtoutfile{\s\s\s  http://www.maths.warwick.ac.uk/agt/AGTVol\thevolumenumber/agt-\thevolumenumber-\thepapernumber.abs.html}
\immediate\write\gtoutfile{\noexpand\\}
\immediate\write\gtoutfile{}
\ifx\theasciiabstract\relax
\immediate\write\gtoutfile{\theabstract}\else
\immediate\write\gtoutfile{\theasciiabstract}\fi
\immediate\write\gtoutfile{}
\immediate\write\gtoutfile{\noexpand\\}
\immediate\write\gtoutfile{}
\immediate\closeout\gtoutfile}}  

\def\maketitlepage{\makeagttitle\makeheadfile}
\let\makeshorttitle\maketitlepage
\let\maketitle\maketitlepage

\def\ifplaintex{\expandafter\ifx\csname documentclass\endcsname\relax}

\def\gtp{{\mathsurround=0pt\it $\cal G\mskip-2mu$eometry \&\ 
$\cal T\!\!$opology $\cal P\!$ublications}}  

\def\recd{{\small Received:\qua\receiveddate\ifx\reviseddate\relax
\else\qquad Revised:\qua\reviseddate\fi\par}} 

\def\lognumber#1{\def\thelognumber{#1}}
\def\volumenumber#1{\def\thevolumenumber{#1}}
\def\volumeyear#1{\def\thevolumeyear{#1}}
\def\papernumber#1{\def\thepapernumber{#1}}
\def\pagenumbers#1#2{\def\startpage{#1}\def\finishpage{#2}}
\def\published#1{\def\publishdate{#1}}

\def\received#1{\def\receiveddate{#1}}

\def\accepted#1{\def\accepteddate{#1}}

\def\covertitle#1{\def\thecovertitle{#1}}

\def\asciiaddress#1{\def\theasciiaddress{#1}}

\long\def\asciiabstract#1{\long\def\theasciiabstract{#1}}

\let\\\par\let\thelognumber\relax\let\thevolumenumber\relax
\let\thepapernumber\relax\let\thevolumeyear\relax\let\startpage\relax
\let\finishpage\relax\let\publishdate\relax\let\receiveddate\relax
\let\reviseddate\relax\let\accepteddate\relax\let\theasciititle\relax
\let\thecovertitle\relax\let\theasciiauthors\relax\let\theasciiaddress\relax
\let\theasciiabstract\relax

\let\theasciiemail\relax

\ifplaintex
\font\logobig=cmssbx10 scaled 3836
\font\logomed=cmssbx10 scaled 2557
\else
\font\logobig=cmssbx10 scaled 4200
\font\logomed=cmssbx10 scaled 2800
\fi

\long\def\makeagttitle{   
\count0=\startpage
\agt\hfill      
\hbox to 45truept{\vbox to 0pt{\vglue -13truept{\logomed A\kern -.37em{\logobig 
T}\kern -.38em G}\vss}\hss}
\break
{\small Volume \thevolumenumber\ (\thevolumeyear)
\startpage--\finishpage\nl
Published: \publishdate}

\vglue .25truein

{\parskip=0pt\leftskip 0pt plus
1fil\def\\{\par\smallskip}{\Large\bf\thetitle}\par\medskip} \vglue
0.05truein

{\parskip=0pt\leftskip 0pt plus 1fil\def\\{\par}{\sc\theauthors}
\par\medskip}%
 
\vglue 0.03truein 

{\small\leftskip 25truept\rightskip 25truept{\bf Abstract}\stdspace\theabstract

{\bf AMS Classification}\stdspace\theprimaryclass
\ifx\thesecondaryclass\relax\else; \thesecondaryclass\fi\par
{\bf Keywords}\stdspace \thekeywords\par}\vglue 7truept

}   

\ifplaintex
\font\phead=cmsl9 scaled 950
\font\pnum=cmbx10 scaled 913
\font\pfoot=cmsl9 scaled 950
\headline{\vbox to 0pt{\vskip -4.5mm\line{\small\phead\ifnum
\count0=\startpage ISSN 1472-2739 (on-line) 1472-2747 (printed)
\hfill {\pnum\folio}\else\ifodd\count0\def\\{ }%
\ifx\theshorttitle\relax\thetitle\else\theshorttitle\fi\hfill{\pnum\folio}
\else\def\\{ and }{\pnum\folio}\hfill\ifx\theshortauthors\relax\theauthors
\else\theshortauthors\fi\fi\fi}\vss}}
\footline{\vbox to 0pt{\vglue 0mm\line{\small\pfoot\ifnum\count0=\startpage
\copyright\ \gtp\hfill\else
\agt, Volume \thevolumenumber\ (\thevolumeyear)\hfill\fi}\vss}}
\else
\font=cmsl9 scaled 1050
\font\lnum=cmbx10 
\font=cmsl9 scaled 1050
\makeatletter
\def\@oddhead{{\small\lhead\ifnum\count0=\startpage ISSN 1472-2739 
(on-line) 1472-2747 (printed)\hfill {\lnum\number\count0}\else\ifodd\count0
\def\\{ }\ifx\theshorttitle\relax \thetitle \else\theshorttitle\fi\hfill
{\lnum\number\count0}\else\def\\{ and }{\lnum\number\count0}
\hfill\ifx\theshortauthors\relax 
\theauthors\else\theshortauthors\fi\fi\fi}}\def\@evenhead{\@oddhead}
\def\@oddfoot{\small\lfoot\ifnum\count0=\startpage\copyright\ \gtp\hfill\else
\agt, Volume \thevolumenumber\ (\thevolumeyear)\hfill\fi}
\def\@evenfoot{\@oddfoot}
\makeatother
\fi
\let\maketitlepage\makeagttitle
\let\makeshorttitle\maketitlepage
\let\maketitle\maketitlepage

\newwrite\gtoutfile
\long\gdef\makeheadfile{  
{\def\\{, }\def\s{ }
\immediate\openout\gtoutfile head.xxx
\immediate\write\gtoutfile{To: math@arxiv.org}
\immediate\write\gtoutfile{Subject: put OR rep NNNNN:ppppp}
\immediate\write\gtoutfile{--text follows this line--}
\immediate\write\gtoutfile{Proxy-for: \ifx\theasciiauthors\relax
\theauthors\else\theasciiauthors\fi\s<\ifx\theasciiemail\relax\theemail\else\theasciiemail\fi>}
\immediate\write\gtoutfile{\noexpand\\}
\immediate\write\gtoutfile{Authors: \ifx\theasciiauthors\relax
\theauthors\else\theasciiauthors\fi}
{\def\\{ }\immediate\write\gtoutfile{Title: \ifx\theasciititle\relax
\thetitle\else\theasciititle\fi}}
\immediate\write\gtoutfile{Subj-class: GT or SG, GR etc}
\immediate\write\gtoutfile{MSC-class: \theprimaryclass\ifx\thesecondaryclass\relax\else, \thesecondaryclass\fi}
\immediate\write\gtoutfile{Journal-ref: Algebr. Geom. Topol. \thevolumenumber\s
(\thevolumeyear) \startpage-\finishpage}
\immediate\write\gtoutfile{Comments: Published by Algebraic and
Geometric Topology at}
\immediate\write\gtoutfile{\s\s\s  http://www.maths.warwick.ac.uk/agt/AGTVol\thevolumenumber/agt-\thevolumenumber-\thepapernumber.abs.html}
\immediate\write\gtoutfile{\noexpand\\}
\immediate\write\gtoutfile{}
\ifx\theasciiabstract\relax
\immediate\write\gtoutfile{\theabstract}\else
\immediate\write\gtoutfile{\theasciiabstract}\fi
\immediate\write\gtoutfile{}
\immediate\write\gtoutfile{\noexpand\\}
\immediate\write\gtoutfile{}
\immediate\closeout\gtoutfile}}  

\def\maketitlepage{\makeagttitle\makeheadfile}
\let\makeshorttitle\maketitlepage
\let\maketitle\maketitlepage

\lognumber{11}
\volumenumber{2}
\volumeyear{2002}
\papernumber{11}
\published{11 April 2002}
\pagenumbers{239}{255}
\received{12 March 2002}
\accepted{22 March 2002}

\usepackage{amsmath}
\usepackage{amssymb}
\usepackage{amscd}

\numberwithin{equation}{section}

\newtheorem{thm}{Theorem}
\newtheorem{thrm}{Theorem}[section]
\newtheorem{prop}[thrm]{Proposition}

\newtheorem{cor}{Corollary}

\theoremstyle{definition}
 \newtheorem*{dfn}{Definition}
 \newtheorem*{qtn}{Question}

\theoremstyle{remark}

\newcommand{\mts}{M\times S^1}
\newcommand{\mtr}{M\times\R}
\newcommand{\mti}{M\times I}
\newcommand{\mtj}{M\times J}
\newcommand{\dt}{{\cdot}}
\newcommand{\rd}{\partial}
\newcommand{\uv}[1]{\rd/\rd #1}
\newcommand{\sst}{\subset}
\newcommand{\pnx}{\pl_n(\xi)}
\newcommand{\ptx}{\tilde{\pl}(\xi)}
\newcommand{\dnx}{\D_n(\xi)}
\newcommand{\dba}[1]{\bar{\D}(#1)}
\newcommand{\dbn}[1]{\bar{\D}_n(#1)}

\newcommand{\dtx}{\tilde{\D}(\xi)}
\newcommand{\tw}[1]{\operatorname{tw}(#1)}
\newcommand{\mtw}[1]{\operatorname{tw_-}(#1)}

\newcommand{\wed}{\wedge}
\newcommand{\inte}{\lrcorner}
\newcommand{\lan}{\langle}
\newcommand{\ran}{\rangle}
\newcommand{\dev}[1]{\Phi_{#1}}
\newcommand{\id}{\operatorname{id}}
\newcommand{\tld}[1]{\tilde{#1}}
\newcommand{\spn}[1]{\operatorname{Span}\lan#1\ran}
\newcommand{\indeng}[2]{\D(\xi,\F_{#1},\F_{#2},n)}

\newcommand{\al}{\alpha}
\newcommand{\be}{\beta}
\newcommand{\vphi}{\varphi}
\newcommand{\eps}{\varepsilon}

\newcommand{\R}{\mathbb{R}}
\newcommand{\D}{\mathcal{D}}
\newcommand{\E}{\mathcal{E}}
\newcommand{\F}{\mathcal{F}}
\newcommand{\pl}{\mathbb{P}}
\newcommand{\cL}{\mathcal{L}}
\newcommand{\N}{\mathbb{N}}
\newcommand{\Z}{\mathbb{Z}}
\newcommand{\nnZ}{\mathbb{Z}_{\ge 0}}

\begin{document}
%
%
\title{Engel structures with trivial characteristic foliations}
\covertitle{Engel structures with trivial\\characteristic foliations}
\author{Jiro Adachi}
\address{Department of Mathematics, 
  Osaka University\\
  Toyonaka Osaka 560-0043, Japan}
\email{adachi@math.sci.osaka-u.ac.jp, jiro@math.upenn.edu}
\secondaddress{Department of Mathematics, 
  University of Pennsylvania\\
  Philadelphia PA 19104, USA}
\asciiaddress{Department of Mathematics, Osaka University\\
  Toyonaka Osaka 560-0043, Japan\\and\\Department of Mathematics, 
  University of Pennsylvania\\
  Philadelphia, PA 19104, USA}

\begin{abstract}
  Engel structures on $\mts$ and $\mti$ are studied in this paper, 
where $M$ is a $3$--dimensional manifold. 
 We suppose that these structures have characteristic line fields parallel 
to the fibres, $S^1$ or $I$. 
 It is proved that they are characterized by contact structures 
on the cross section $M$, the twisting numbers, and Legendrian foliations 
on both ends $M\times \rd I$ in the case of $\mti$. 
\end{abstract} 

\asciiabstract{Engel structures on M x S^1 and M x I are studied in
this paper, where M is a 3-dimensional manifold.  We suppose that these
structures have characteristic line fields parallel to the fibres,
S^1 or I.  It is proved that they are characterized by contact
structures on the cross section M, the twisting numbers, and
Legendrian foliations on both ends M x dI in the case of
M x I.}

\maketitle

\section{Introduction}
  In any category of geometry, 
exotic structures are interesting objects to study. 
 In addition, the classification is an important problem. 
 Engel structures are largely dominated by the characteristic foliations. 
 This paper is devoted to the characterization of Engel structures 
on $\mts$ and $\mti$, where $M$ is a $3$--dimensional manifold, 
with the same characteristic foliation as the standard Engel structure. 

  A maximally non-integrable distribution of rank $2$ 
on a $4$--dimensional manifold is called an \emph{Engel structure\/} 
(see Section~\ref{sec:dfns} for the precise definition). 
 Engel structures have have no local invariant, 
similarly to contact and symplectic structures. 
 There exists a local normal form, 
written as a kernel of two differential $1$--forms, 
\begin{equation}  \label{eq:normal}
  dy-z\dt dx=0, \qquad dz-w\dt dx=0, 
\end{equation}
where $(x,y,z,w)$ is a coordinate system. 
 Further, it is known that this property occurs only 
on line fields, contact structures, even-contact structures, 
and Engel structures, among regular tangent distributions on manifolds 
(see \cite{mont1}, \cite{vege}). 
 An \emph{even-contact structure\/} is a hyperplane field 
on an even dimensional manifold with maximal non-integrability. 
 The fact above indicates the importance of the study of Engel structures. 
 However, different from contact structures, 
global stability does not hold for Engel structures, 
that is, the Gray type theorem does not hold. 
 There exists a line field $\cL(\D^2)$, which is canonically defined 
by a given Engel structure $\D$, 
especially by an even-contact structure $\D^2:=\D+[\D,\D]$ 
(see Section~\ref{sec:dfns}). 
 It is called the \emph{characteristic line field\/} of $\D$ or $\D^2$. 
 The foliation $L(\D^2)$,  
whose leaves are integral curves of the characteristic line fields, 
is called the \emph{characteristic foliation} of $\D$ or $\D^2$. 
 This is the obstruction to global stability. 
 It is proved by A.~Golubev and R.~Montgomery in \cite{gol} and \cite{mont2} 
that a deformation of an Engel structure is realized by an isotopy 
if it fix the characteristic line field. 
 If the above two differential $1$--forms~\eqref{eq:normal} 
are defined globally on $\R^4$, 
the obtained Engel structure is called the \emph{standard\/} Engel structure 
on $\R^4$. 
 Let $\D_{st}$ denote it. 
 The standard Engel structure has its characteristic line field 
spanned by a vector field in the $w$ direction, 
$\cL({\D_{st}}^2)=\spn{\uv{w}}$. 
 V.~Gershkovich constructs in \cite{ger} examples of Exotic Engel structures 
on $\R^4$ with periodic orbits in their characteristic line fields. 
 To study exotic Engel structures 
with the same characteristic line field as the standard Engel structure 
is a motivation for this paper. 
 
  In order to state the results of this paper, we prepare some basic notions. 
 Let $\xi$ be a contact structure on a $3$--dimensional manifold $M$. 
 Let us construct a new $4$--dimensional manifold $\pl(\xi)$ from $\xi$. 
 A contact structure $\xi$ is a certain $2$--plane field. 
 Then the manifold $\pl(\xi)$ is obtained 
by the fibrewise projectivization of $\xi$: 
$\pl(\xi)=\cup_{p\in M} \pl(\xi_p)$, 
where $\pl(\xi_p)$ is the projectivization of the tangent plane $\xi_p$. 
 This manifold has a structure of the $S^1$--bundle over $M$, 
because each fibre is equivalent to $S^1$: $\pl(\xi_p)\cong \R P^1\cong S^1$. 
 We can regard a point $l\in \pl(\xi_p)$ of a fibre 
as a line on the plane $\xi_p$ through the origin. 
 A $2$--plane $\D(\xi)_{(p,l)}$ at a point $(p,l)\in\pl(\xi)$ is induced 
by the pull back of a line $l\sst \pl(\xi_p)$ 
by the projection $\pi\co \pl(\xi)\to M$. 
 Then we obtain a $2$--plane field $\D(\xi)$ on $\pl(\xi)$. 
 It is known that this $2$--plane field $\D(\xi)$ is an Engel structure 
on $\pl(\xi)$ (see \cite{mont2}). 
This Engel manifold $(\pl(\xi),\D(\xi))$ is called a \emph{prolongation\/} 
of a contact $3$--manifold $(M,\xi)$ 
(see Section~\ref{sec:prolong} for the definition). 
 Furthermore, we consider its fibrewise $n$--fold covering 
and the corresponding Engel structure. 
 Let $(\pl_n(\xi), \D_n(\xi))$ denote them. 
 We regard $(\pl(\xi), \D(\xi))$ as $(\pl_1(\xi), \D_1(\xi))$, if necessary. 
 We note that the prolonged manifold $\pl_n(\xi)$ is diffeomorphic to $\mts$ 
if the given contact structure $\xi$ is trivial as plane field. 
 There is a canonical identification $\psi_n\co \pl_n(\xi)\to \mts$ 
in this case. 
 According to this identification, we obtain an Engel structure 
$\dbn{\xi}:=(\psi_n)_{\ast}\D_n(\xi)$ on $\mts$ from a contact structure $\xi$ 
on $M$ and a natural number $n\in\N$. 
 Prolongations and deprolongations are discussed in Section~\ref{sec:prolong}. 

  In this paper, Engel structures on $\mts$ and $\mti$ 
with trivial characteristic foliations are investigated. 
 We say a characteristic foliation $L$ is \emph{trivial} 
if it is isotopic to a foliation which consists of leaves 
$\{pt\}\times S^1\sst\mts$ or $\{pt\}\times I\sst\mti$. 
 In the following, we suppose that these isotopies have been applied. 
 Namely, we assume in the following 
that a trivial characteristic foliation consists 
of leaves $\{pt\}\times S^1$ or $\{pt\}\times I$. 
 In this case, we can define invariants for Engel structures. 
 Let $\D$ be an Engel structure on $\mts$ 
with a trivial characteristic foliation. 
 An Engel structure, as a $2$--plane field, is spanned 
by the characteristic line field and a line field twisting along leaves 
of the characteristic foliation (see Section~\ref{sec:dfns}). 
 Now each leaf of characteristic foliation is a fibre 
$\{pt\}\times S^1 \sst \mts$. 
 Then we can define an invariant, 
the \emph{twisting number\/} $\tw{\D}$ of $\D$, 
by counting the rotation along a leaf with respect to a fixed framing of $\xi$.
 Similarly we, define the \emph{minimal twisting number\/} $\mtw{\D}$ 
for Engel structures on $\mti$. 
(see Section~\ref{sec:devmap} for the definitions.)

  Next, we consider the induced Legendrian foliations. 
 Let $\D$ be an Engel structure on $\mts$ or $\mti$ 
with a trivial characteristic foliation $L(\D^2)$. 
 Let us identify cross sections of $\mts$ and $\mti$ with $M$ itself
by the standard projection. 
 This projection is the projection along the characteristic foliation $L(\D^2)$
too. 
 We note that this $M$ is transverse to the characteristic foliation $L(\D^2)$.
 It is known that the even-contact structure $\D^2$ 
induce a contact structure $\xi(\D^2)$ on $M$ 
(see \cite{mont2} and \cite{ger}). 
 It does not depend on the choice of the cross section 
(see Section~\ref{sec:prolong}).
 Similarly the Engel distribution $\D$ induce an $1$--dimensional foliation 
$\F(M,\D)$ on $M$ by the integral of the line field 
defined by the intersection of $\D$ and $M$. 
 We note that the leaves of this foliation $\F(M,\D)$ are tangent 
to the induced contact structure $\xi(\D^2)$ everywhere. 
 Curves in contact $3$--manifolds, which are tangent to the contact structures 
everywhere, are called \emph{Legendrian curves}. 
 We call this foliation $\F(M,\D)$ the \emph{induced Legendrian foliation}. 

  Now, we are ready to state the results of this paper. 
 First, we consider Engel structures on $\mts$ 
with trivial characteristic foliations. 

%
%
\begin{thm} \label{thms}
\emph{(1)}\qua Let $\xi$, $\zeta$ be parallelizable contact structures, 
that is, they have global framings. 
 If the prolonged Engel structures $\bar{\D}(\xi)$ and $\bar{\D}(\zeta)$ 
are isotopic preserving the characteristic foliation, 
then the contact structures $\xi$ and $\zeta$ are isotopic.

\emph{(2)}\qua For any Engel structure $\D$ on a $4$--dimensional manifold $\mts$ 
with a trivial characteristic foliation, 
there exist a contact structure $\xi$ on $M$ and a natural number $n\in\N$, 
for which the $n$--fold prolonged Engel structure 
$\dbn{\xi}$ of\/ $\xi$ on $\mts$ is isotopic to $\D$. 
\end{thm}
 In other words, this theorem implies the following. 
%
%
\begin{cor} 
  Engel structures on $\mts$ with trivial characteristic foliations 
are characterized, up to isotopies, 
by isotopy classes of contact structures on $M$ 
and the twisting number $\tw{\D}\in \N$. 
\end{cor} 
 We note that we can construct an Engel structure $\dbn{\xi}$ on $\mts$ 
for any isotopy class $[\xi]$ of of contact structures on $M$ 
which are trivial as plane fields, 
and a natural number $n\in \N$. 

  Next, we show that Engel structures on $\mti$ 
with the trivial characteristic foliation are determined 
by the induced contact structures, the minimal twisting numbers, 
and the induced Legendrian foliations on both ends $M\times\rd I$. 
%
%
\begin{thm} \label{thmi}
\emph{(1)}\qua Let $\D$ and $\tilde{\D}$ be Engel structures on $\mti$ 
with the trivial characteristic foliations. 
 If they have the same induced contact structure, 
induced Legendrian foliations on both ends $M\times\rd I$, 
and minimal twisting number, 
then they are isotopic relative to the ends.

\emph{(2)}\qua Let $\xi$ be a parallelizable contact structure 
on a $3$--manifold $M$, 
$(\F_0$, $\F_1)$ a pair of Legendrian foliations on $(M,\xi)$, 
and $n\in\nnZ$ a non-negative integer. 
 Then there exists an Engel structure $\D=\D(\xi,\F_0,\F_1,n)$ on $\mti$, 
which has the induced contact structure $\xi(\D)=\xi$, 
the induced Legendrian foliations $\F(M\times\{i\},\D)=\F_i, i=0,1$, 
and the minimal twisting number $\mtw\D=n$. 
\end{thm} 
 This theorem implies the following. 
%
%
\begin{cor}
  Engel structures on $\mti$ with the trivial characteristic foliation 
are determined by the induced contact structures 
and the induced Legendrian foliations on both ends $M\times\rd I$ and 
the minimal twisting numbers. 
\end{cor} 

\medskip

  This paper was established while the author was visiting 
at Stanford University and the University of Pennsylvania. 
 He would like to express his gratitude 
to Professor Yakov Eliashberg and Professor John Etnyre 
for some discussions and the hospitality. 
 This work was partially supported 
by JSPS Research Fellowship for Young Scientists. 

\section{Preliminary}
\subsection{Definitions and important properties.} \label{sec:dfns}
  Let $W$ be a $4$--dimensional manifold. 
 A \emph{distribution\/} of rank $2$ or a $2$--\emph{plane field\/} $\D$ on $W$
is a distribution of $2$--dimensional tangent planes $\D_p\sst T_pW$, $p\in W$.
 It is considered as a rank $2$ subbundle of the tangent bundle $TW$. 
 We can think of $\D$ as a locally free sheaf of smooth vector fields on $M$. 
 Let $[\D,\D]$ denote the sheaf generated by all Lie brackets $[X,Y]$ 
of vector fields $X, Y$ on $M$, which are cross sections of $\D$. 
 We set $\D^2:= \D+[\D,\D]$ and $\D^3:= \D^2+[\D,\D^2]$. 
 The Engel structure is defined as follows. 
%
%
\begin{dfn}
  A distribution $\D$ of rank $2$ on a $4$--dimensional manifold $W$ is called 
an \emph{Engel structure\/} if it satisfies, 
\begin{equation}  \label{eq:econdi}
  \dim\D^2 = 3, \qquad \dim\D^3 = 4, 
\end{equation}
at any point $p\in W$. 
\end{dfn} 
 We note that $\D^2$ is a distribution of rank $3$ 
and $\D^3$ is the tangent bundle $TW$ itself, if $\D$ is an Engel structure. 
 This corank $1$ distribution $\D^2$ is an even-contact structure on $W$. 
 Let $\E$ denote it. 
 An \emph{even-contact structure} is, by definition, 
a corank $1$ distribution $\E$ on an even-dimensional manifold, 
which is defined, at least locally, by $1$--form $\theta$ with a property 
that $ \theta \wed (d\theta)^{n/2-1}$ is never-vanishing $(n-1)$--form, 
where $n$ is the dimension of the manifold. 
 An even-contact structure $\E$ on a $4$--manifold $W$ has a characteristic 
$1$--dimension. 
 We define a rank $1$ subdistribution $\cL(\E)$ of $\E$ 
by $[\cL(\E),\E]\sst \E$. 
 In other words, it is a maximal integrable subdistribution of $\E$. 
 In this case, its rank is $1$. 
 It is called the \emph{characteristic line field\/} 
of $\E=\D^2$ or sometimes of $\D$. 
 We call the $1$--dimensional foliation 
obtained by integrating the characteristic line field 
the \emph{characteristic foliation\/} of $\D^2$ or $\D$. 
 We note that an Engel structure $\D$ 
should include its characteristic line field 
to satisfy the Engel condition~\eqref{eq:econdi}. 

  Contact structures are sometimes described in terms of contact forms. 
 Similarly, Engel structures are described 
in terms of pairs of differential $1$--forms.  
 A pair of $1$--forms $(\al, \be)$ on a $4$--manifold $W$ 
is called an \emph{Engel pair of\/ $1$--forms\/} 
if it satisfies the following three conditions, 
\begin{enumerate}
 \item[(1)] $\al\wed\be\wed d\al$ never vanishes, 
 \item[(2)] $\al\wed\be\wed d\be \equiv 0$, 
 \item[(3)] $\be\wed d\be$ is a never-vanishing $3$--form. 
\end{enumerate}
 It is known that the distribution $\D:=\{\al=0,\ \be=0\}$ defined 
by an Engel pair of $1$--forms is an Engel structure (see \cite{ger}). 
 Under these conditions above, 
the $1$--form $\be $ defines an even-contact distribution as $\D^2=\{\be =0\}$,
and the characteristic line field $\cL(\D^2)$ is defined 
as a kernel of the $3$--form $\be\wed d\be$. 
 For an even-contact distribution $\E=\{\be=0\}$, 
a vector field $X_0$ is called the \emph{characteristic vector field\/} 
of $\E$, 
if it satisfies $X_0\inte \Omega =\be\wed d\be$ for some volume form $\Omega$ 
on $W$. 
 We note that it generates the characteristic line field and foliation of $\E$.

  At the end of this section, we introduce important properties 
of Engel structures. 
 They play important roles in the proof of Theorems. 
 The following theorem is proved by R.~Montgomery in \cite{mont2}, 
after E.~Cartan's theory. 
%
%
\begin{thrm} \label{thm:mont}
  Let $\D_0,\ \D_1$be Engel structures 
defined near an embedded $3$--manifold $M$ in a $4$--manifold $W$. 
 We assume that their characteristic foliation $L(\D_i^2)$, $i=0,1$, 
are transverse to $M$, 
and that these Engel structures $\D_i$, $i=0,1$, 
have the same induced Legendrian foliation on $M$: 
$\F(M,\D_0)=\F(M,\D_1)$. 
 Then they are locally isotopic along $M\sst W$. 
\end{thrm}
 A slight generalization of this theorem is shown in \cite{adachi}. 
 The case where the embedded $3$--manifold is not transverse 
to the characteristic foliation is considered in \cite{adachi}. 
 The following theorem is proved in \cite{mont2} and \cite{gol}. 
 It is well known that contact structures on a closed orientable $3$--manifold 
are isotopic if they are homotopic among contact structures (the Gray Theorem).
 The following theorem is the Gray type theorem for Engel structures. 
%
%
\begin{thrm} \label{thm:gol}
  Let $\D_t$, $t\in[0,1]$, be a family of Engel structures 
on a closed orientable $4$--manifold $W$. 
 We assume that it has the fixed characteristic foliation: $L(\D_t^2)\equiv L$.
 Then there exists a family $\phi_t\co W\to W$, $t\in[0,1]$, 
of global diffeomorphisms which satisfies 
$(\phi_t)_{\ast}\D_t=\D_0$, $(\phi_t)_{\ast}L=L$. 
\end{thrm}

\subsection{Prolongation procedures of contact manifolds.} \label{sec:prolong}
  The notion of prolongation is introduced by E.~Cartan 
in the theory of exterior differential systems (see \cite{cartan}, \cite{bcg}).
 We consider the prolongations of contact structures on $3$--manifolds. 
 Let $\xi$ be a contact structure on a $3$--manifold $M$, 
namely a certain $2$--plane field. 
 We construct a new $4$--dimensional manifold from $\xi$ 
by fibrewise projectivizations, 
\begin{equation*}
  \pl(\xi):= \bigcup_{p\in M} \pl(\xi_p), 
\end{equation*}
where $\pl(\xi_p)$ is a projectivization of a tangent plane $\xi_p$. 
 A point of $\pl(\xi)$ can be regarded as a line $l$ 
in the contact plane $\xi_p$ through the origin. 
 The constructed $4$--manifold $\pl(\xi)$ has a structure of a circle bundle 
over $M$. 
 Let $\pi\co \pl(\xi)\to M$ be its projection. 
 This $4$--manifold is endowed with a $2$--plane field $\D(\xi)$ 
induced naturally as follows. 
 We define $2$--plane $\D(\xi)_q\sst T_q(\pl(\xi))$ at $q=(p,l)\in\pl(\xi)$. 
 A point $q=(p,l)\in\pl(\xi)$ is regarded as a pair of a point $p\in M$ and 
a tangent line $l\sst \xi_p\sst T_pM$. 
 Then we set $\D(\xi)_q:= (d\pi^{-1})_ql$. 
 Thus we obtain a $2$--plane field $\D(\xi)$ on a $4$--manifold $\pl(\xi)$. 
 We call this $(\pl(\xi), \D(\xi))$ the \emph{prolongation\/} 
of a contact structure $\xi$ on a $3$--manifold $M$. 
 It is known that the prolongation $(\pl(\xi), \D(\xi))$ is an Engel manifold 
(see \cite{mont2}). 
 We note that the prolonged manifold $\pl(\xi)$ is diffeomorphic to $\mts$ 
if the contact structure $\xi$ belongs to the trivial class as plane fields. 
 In this paper we consider this case especially. 

  Furthermore, we consider some variants of prolongations. 
 Let $(\pl(\xi), \D(\xi))$ be the prolongation of a contact structure $\xi$ on 
a $3$--manifold $M$. 
 We define a new $4$--manifold $\pl_n(\xi)$ 
by fibrewise $n$--fold covering of $\pl(\xi)$, 
\begin{equation*}
  \pl_n(\xi):= \bigcup_{p\in M} \pl_n(\xi_p), 
\end{equation*}
where $\pl_n(\xi_p)$ is an $n$--fold covering space of $\pl(\xi_p)$. 
 Let $\vphi_n\co \pl_n(\xi)\to \pl(\xi)$ 
be a fibrewise covering bundle mapping. 
 We obtain a corresponding Engel structure $\D_n(\xi)$ on $\pl_n(\xi)$: 
$(\vphi_n)_{\ast}\D_n(\xi)=\D(\xi)$. 
 This pair $(\pl_n(\xi), \D_n(\xi))$ is called 
the \emph{$n$--fold prolongation\/} of a contact structure 
on a $3$--manifold $M$. 
 According to this notation, we have 
$(\pl(\xi),\D(\xi))=(\pl_1(\xi),\D_1(\xi))$. 
 When the given contact structure $\xi$ on a $3$--manifold $M$ 
is trivial as plane fields, 
$\pl_n(\xi)$ is diffeomorphic to $\mts$ for any $n\in \N$. 
 Let $V_0, V_1$ be vector fields on $M$ spanning $\xi$: 
$\xi=\spn{V_0, V_1}$. 
 There is a canonical identification $\psi_n\co \mts\to \pnx$ 
defined as follows, 
\begin{equation*}
  \psi_n\co (p,\theta)\mapsto 
          \left(p, l=\left[(V_0)_p\dt\cos\left(\frac{n\theta}{2}\right)
                  +(V_1)_p\dt\sin\left(\frac{n\theta}{2}\right)\right]\right), 
\end{equation*}
where $l\in\pl_n(\xi_p)$ is a line in $\xi_p$ defined by a vector field 
$(V_0)_p\dt\cos(n\theta/2)+(V_1)_p\dt\sin(n\theta/2)$. 
 Then we obtain a corresponding Engel structure on $\mts$ by setting 
$\dbn{\xi}:=(\psi_n^{-1})_{\ast}\D_n(\xi)$, $n=1,2,3,\ldots$. 
 In addition, we consider a fibrewise universal covering 
of a prolongation $\pl(\xi)$, 
and let $(\tilde{\pl}(\xi),\tilde{\D}(\xi))$ denote it. 

  Next, we consider the deprolongation procedure, 
the inverse of the prolongation, in a sense. 
 Similarly to the above, we consider deprolongation of Engel structures 
especially. 
 Let $\D$ be an Engel structure on a $4$--manifold $W$, 
$\E:=\D^2$ its even-contact structure, 
and $L(\E)$ its characteristic foliation. 
 We consider the leaf space $W/L(\E)$ 
and its projection $\pi\co W\to W/L(\E)$. 
 The foliation $L(\E)$ is said to be \emph{nice}, according to \cite{mont2},  
if $W/L(\E)$ is a smooth $3$--manifold and $\pi$ is a submersion. 
 We suppose here that $L(\E)$ is nice. 
 We set $\xi(\E):=\pi_{\ast}\E$. 
 It is a $2$--plane field on $W/L(\E)$, 
which is well defined because the characteristic vector field $X_0$, 
along $L(\E)$, preserves the even-contact structure $\E$. 
 In fact, the even-contact structure $\E=\D^2$ is determined 
by the second $1$--form $\beta$ 
of the Engel pair of $1$--forms $(\alpha,\beta)$. 
 Since the characteristic vector field $X_0$ is defined 
so that $X_0\inte(\beta\wedge d\beta)=0$ (see Section~\ref{sec:dfns}), 
we have $L_{X_0}\beta=f\dt\beta$, for some function $f$. 
 This implies that $X_0$ preserves $\E$. 
 Therefore, we have $\xi(\E)_{\pi(p)}:=(d\pi)_p(\E_p)=(d\pi)_q(\E_q)$ 
for any point $q$ on the same leaf of $L(\E)$ as $p$ 
because $\pi$ is the projection along $L(\E)$ or $X_0$. 
 It is known that this distribution $\xi=\xi(\E)$ is a contact structure 
on $W/L(\E)$ (see \cite{mont2}, \cite{ger}). 
 We call this $(W/L(\E),\xi(\E))$ the \emph{deprolongation\/} 
of the Engel structure $\D$. 

  Let $(W,\D)$ be an Engel manifold 
with the characteristic foliation $L(\D^2)$,
and $M\sst (W,\D)$ an embedded $3$--manifold. 
 We assume that $M$ is transverse to the characteristic foliation $L$. 
 Then we can take a neighborhood $U$ of $M$ as a flow-box for $L$. 
 In this neighborhood $U$, the foliation $L$ is nice in the sense above. 
 Thus we can apply the deprolongation procedure. 
 In this case, we can identify the leaf space $U/L(\D^2)$ with $M$. 
 Then the obtained contact structure is $\pi_{\ast}\D^2=TM\cap \D^2$. 
 In the case where the Engel manifold is a prolongation $(\pl(\xi),\D(\xi))$ 
of a contact structure $\xi$ on a $3$--manifold $M$ and embedded $3$--manifold 
is a cross section $M_{\theta}$, 
the characteristic foliation $L(\D(\xi)^2)$ is nice globally. 
 Then, the leaf space $\pl(\xi)/L(\D^2)$ is identified with $M_{\theta}=M$ 
and the obtained structure is the original $\xi$. 

\subsection{Twisting property and Development mappings.} \label{sec:devmap}
  The development mapping is a local Engel diffeomorphism or an immersion 
into a prolongation $(\pl(\xi),\D(\xi))$ or $(\ptx,\dtx)$, 
introduced in \cite{mont2}. 
 It is constructed by a property that an Engel structure is twisting 
along leaves of its characteristic foliation. 
 First, we observe the twisting condition of Engel structures. 
 Let us recall that an Engel structure $\D$ contains 
its characteristic line field $\cL=\cL(\D^2)\sst \D$ of $\D^2=:\E$. 
 Let $\tilde{\D}$ be another rank $2$ distribution which is contained in $\E$ 
and contains $\cL=\cL(\E)$: $\cL\sst \tilde{\D}\sst \E$. 
 Then the \emph{twisting condition}: 
\begin{equation} \label{eq:twcondi}
   \tilde{\D}+ [\cL,\tilde{\D}]=\E
\end{equation}
implies that $\tilde{\D}$ is Engel, 
that is, the Engel condition \eqref{eq:econdi}. 
 Let $X_0$ be a characteristic vector field of $\E$, and $V$ a vector field 
which forms a basis of $\tilde{\D}$ with $X_0$. 
 We note $[V,\E]\not\subset \E$ 
because $V$ is independent of the integrable subdistribution $\cL(\E)\sst \E$. 
 Then the derived distributions are 
\begin{align*}
  \tilde{\D}^2 &=\tilde{\D}+[\cL,\tilde{\D}]+[V, \tilde{\D}] 
                =\E, \\ 
  \tilde{\D}^3 &=\E+[\cL,\E]+[V,\E] 
                =\E+[V,\E], 
\end{align*} 
where we use the fact that the rank of $\tilde{\D}$ is $2$. 
 Therefore we obtain $\dim \tilde{\D}^2=3$, $\dim \tilde{\D}^2=4$, 
namely, the Engel condition \eqref{eq:econdi}.  

  Using this twisting condition, we construct a mapping 
from a domain in an Engel manifold with a nice characteristic foliation
to a prolongation of some contact structure. 
 Let $(W,\D)$ be an Engel manifold, 
and $L=L(\D^2)$ the characteristic foliation. 
 Let $U\sst (W,\D)$ be a domain where the characteristic foliation $L$ is nice 
in the sense above, and $\pi\co U\to U/L$ the projection to the leaf space. 
 We set $\xi:=\pi_{\ast}(\D^2)$. 
 It is a deprolonged contact structure on $U/L$ from $\D$. 
 Let $l\sst W$ be a leaf of $L$, and $q\in l\cap U$ a point. 
 From another point of view, $l$ is a point of the leaf space $U/L$. 
 For a point $q\in l$, there corresponds a tangent $2$--plane $\D_q\sst T_qU$. 
 Since an Engel distribution $\D$ contains the characteristic line field 
$\cL=\cL(\D^2)$, 
and is contained in $\D^2$: $\cL\sst \D\sst \D^2$, 
there corresponds a tangent line $d\pi_q(\D_q)\sst d\pi_q(\D^2_q)=\xi_l$, 
that is, a point of $\pl(\xi_l)$. 
 In this way, we obtain a mapping from $l$ to $\pl(\xi_l)$. 
 As the domain $U$ with a nice foliation is regarded as a union 
of leaves $\cup_{l\cap U\ne\emptyset}(l\cap U)$, 
we obtain a map from a nice domain $U$ to a prolongation $\pl(\xi)$. 
 We call this mapping $\dev{\D}\co U\to \pl(\xi)$ 
the \emph{development mapping\/} associate to the Engel structure $\D$. 
 The twisting condition \eqref{eq:twcondi},  and the argument following it, 
ensure that this development mapping is diffeomorphic locally. 
 From the construction, the development mapping $\dev{\D}$ preserves 
the characteristic line field and another line field in the Engel distribution 
twisting in the sense of the twisting condition, 
with respect to the given $\D$ on $U$ 
and the prolongation $\D(\xi)$ on $\pl(\xi)$. 
 Therefore the development mapping is an Engel diffeomorphic locally 
(see \cite{mont2}). 
 We consider the lift $\tilde{\Phi}_{\D}\co U\to\tld{\pl}(\xi)$ 
of the development mapping. 
 We also call this the development mapping. 

  Now, we define the twisting number for Engel structures on $\mts$ and $\mti$,
with trivial characteristic foliations. 
 We begin with $\mts$. 
 Let $\D$ be an Engel structure on $\mts$ 
with a trivial characteristic foliation, 
and $\xi=\xi(\D)$ its deprolongation. 
 Let $l\in L(\D^2)$ be a leaf of the characteristic foliation 
corresponding to a point $p\in M=(\mti)/L(\D^2)$. 
 We note that $l$ is diffeomorphic to $S^1$. 
 Then we obtain a mapping from $l\cong S^1$ 
to $\pl(\xi_p)\cong\R P^1\cong S^1$, defined as 
$\theta\mapsto\D_{(p,\theta)}\cap T_pM$, 
which can be regarded as a mapping $S^1\to S^1$. 
 The \emph{twisting number\/} $\tw\D$ of $\D$ is defined 
as the degree of this mapping. 
 We suppose that the orientation of a fibre $l$ is defined 
by the characteristic vector field $X_0$. 
 Considering the basis $(v_0,v_1)$ of $\xi$ such that $[X_0,v_0]=v_1$, 
we obtain an orientation of $\xi$, 
and then that of $S^1\cong\pl(\xi_p)$. 
 We consider the above degree with respect to those orientations. 
 We note that it is independent of the choice of points $p\in M$. 
 In other words, the Engel structure $\D$ with the twisting number $\tw\D=n$ 
is Engel diffeomorphic to $(\pnx,\dnx)$ by the development mapping, 
where $\xi=\xi(\D)$. 
 Next, we define the minimal twisting number for Engel structures on $\mti$ 
with trivial characteristic foliations. 
 Let $\bar{\D}$ be an Engel structure on $\mti$ 
with the trivial characteristic foliation, 
and $\bar{\xi}=\xi(\bar{\D})$ its deprolongation. 
 We take a pair $(V_0,V_1)$ of non-vanishing vector fields on $M$, 
which is a positively oriented basis of $\xi$ as above, 
such that $V_0$ defines the induced Legendrian foliation 
$\F_0=\F(M\times\{0\},\bar{\D})$. 
 Then we obtain a family $g_t\co M\to\R$ of functions, which satisfies 
that $g_0\equiv 0$ 
and that the vector fields $(V_0\dt\cos(g_t\pi)+V_1\dt\sin(g_t\pi))$ 
define the line fields $\bar{\D}\cap TM_t$, by identifying $M_t$ with $M$. 
 We call the non-negative integer $n\in\nnZ$ 
such that $n\le\min_{p\in M}g_1<n+1$ the \emph{minimal twisting number\/} 
of $\bar{\D}$. 
 Let $\mtw{\bar{\D}}$ denote it. 
 We note that it is independent of the choice of $V_1$ 
and the orientation of $V_0$.

\section{Characterization on $M\times S^1$} \label{sec:pr1}
  In this section, we prove Theorem~\ref{thms}. 
 First, we show $(1)$ of the theorem. 

\proof[Proof of Theorem \ref{thms}-(1)]
  Let $\vphi_t\co \mts\to \mts$, $t\in [0,1]$, be a given isotopy 
between $\dba{\xi}$ and $\dba{\zeta}$, 
which satisfies $\vphi_0=\id$, $\vphi_{1\ast}(\dba{\zeta})=\dba{\xi}$, 
and $(\vphi_t)_{\ast}L=L$, 
where $L=L(\dba{\xi})=L(\dba{\zeta})$ is a trivial characteristic foliation. 
 We set $\D_t:=(\vphi_t)_{\ast}\dba{\zeta}$. 
 It is a family of Engel structures with the same characteristic foliation $L$.
 Let $\E_t$ denote an even-contact structure $\D_t^2$ on $\mts$. 
 We set $M(t):=\vphi_t(M_0)\sst\mts$, 
where $M_0=M\times\{0\}\sst \mts$ is the zero section. 
 We note that $M(t)$ intersects every leaf of $L$ transversely, 
for any $t\in[0,1]$. 
 As mentioned in Section~\ref{sec:prolong}, 
$\xi_t:=TM(t)\cap\E_t$ is a contact structure on $M(t)$. 
 Each $M(t)$ is sent to $M_0$ by a time one mapping of a vector field 
$f_t\dt (\uv{\theta})=:W_t$ for some function $f_t$, 
where $\theta\in S^1$ is the coordinate of a fibre. 
 Let $\psi^t_1\co \mts\to \mts$ be this family of time one mappings. 
 Each mapping $\psi^t_1$ preserves each even-contact structure $\E_t$, 
because $\uv{\theta}$ is the characteristic vector field for each $\E_t$ 
(See Section~\ref{sec:prolong}). 
 Then, the restriction of $\psi^t_1\circ\vphi_t\co\mts\to\mts$ to $M_0=M$ 
is the required isotopy. 
 In fact, we have $\psi^t_1\circ \vphi_t(M_0)=M_0$, 
$\psi^0_1\circ \vphi_0=\id$, and 
\begin{align*}
   (\psi^1_1\circ \vphi_1)_{\ast}\zeta 
 &={\psi^1_1}_{\ast}\circ{\vphi_1}_{\ast}(TM_0\cap\dba{\zeta}^2) \\
 &={\psi^1_1}_{\ast}(TM(1)\cap\dba{\xi}^2) \\
 &=TM_0\cap\dba{\xi}^2=\xi.\tag*{\qed} 
\end{align*}

  Next, we show $(2)$ of Theorem~\ref{thms}.

\proof[Proof of Theorem \ref{thms}-(2)] 
  Let $\D$ be a given Engel structure on $\mts$ 
with the trivial characteristic foliation $L$. 
 We set $\E:=\D^2$. 
 It is an even-contact structure on $\mts$. 
 We identify the leaf space $(\mts)/L$ with $M_0=M\times \{0\}\sst \mts$. 
 Then the projection $\mts\to (\mts)/L$ is the first projection 
$\pi\co \mts\to M=M_0$. 
 In terms of Section~\ref{sec:prolong}, this foliation $L$ is nice. 
 According to a deprolongation procedure, 
we obtain a contact structure $\xi:=\pi_{\ast}\E$ on $M$ 
(see Section~\ref{sec:prolong}). 
 Let $n\in\N$ be the twisting number of $\D$: $\tw{\D}=n$. 
 We show that the given structure $\D$ is isotopic,  
preserving the characteristic foliation $L$, 
to an $n$--fold prolongation $\dbn{\xi}$ of a deprolongation $\xi$. 
 Now, we consider a development mapping 
$\dev{\D}\co (\mts,\D)\to (\pl(\xi),\D(\xi))$, 
defined in Section~\ref{sec:devmap}. 
 In this case, the characteristic foliation is nice globally, 
and each leaf is $S^1$. 
 Therefore, it is immersed into a prolongation $(\pl(\xi),\D(\xi))$ 
covering $n$ times. 
 In fact, each leaf $l\in L$ is immersed 
onto each leaf of the characteristic foliation of $\D(\xi)^2$, 
winding $n$ times, 
because a line $d\pi(\D)\sst \xi$ twists $n/2$ times along $l$. 
 Then we obtain an Engel diffeomorphism 
$\tilde{\Phi}(\D)\co (\mts,\D)\to (\pnx,\dnx)$ 
as a lift for a fibrewise covering mapping 
$\vphi_n\co (\pnx,\dnx)\to (\pl(\xi),\D(\xi))$. 
 By a canonical identification $\psi_n$, we obtain an Engel diffeomorphism 
$\psi_n^{-1}\circ\tilde{\Phi}(\D)\co (\mts,\D)\to (\mts,\dbn{\xi})$. 
 We note that it is isotopic to identity because it is a bundle map 
over identity between product bundles. 
\endproof

  This Theorem~\ref{thms} implies that Engel structures on $\mts$
with the trivial characteristic foliation are classified, up to isotopy, 
if contact structures on the $3$--manifold $M$ are classified. 
 The classification of contact structure has been an important problem 
for a long time. 
 There are some results on this subject 
(for example, \cite{classot}, \cite{20years}, \cite{kanda}, \cite{etnyre}, 
\cite{giroux}, \cite{honda}). 
 Here we take $S^3$ for example. 

\noindent \textbf{Example}\qua
  Let us consider Engel structures on $S^3\times S^1$ 
with trivial characteristic foliations. 
 Contact structures on $S^3$ are classified, up to isotopy, 
to the following structures (see \cite{classot}, \cite{20years}): 
a tight structure $\zeta$, overtwisted structures $\xi_m$, $m\in\Z$. 
 Then, any Engel structure on $S^3\times S^1$ 
with a trivial characteristic foliation is isotopic to one of the followings: 
\begin{equation*}
  \dbn{\zeta},\quad \dbn{\xi_m},\quad n\in\N,\ m\in\Z. 
\end{equation*}
 These Engel structures are not isotopic each other.

\section{Characterization on $\mti$}
  In this section, we prove Theorem~\ref{thmi}. 
 First, we show that the given pair of Legendrian foliations 
on the contact manifold $(M,\xi)$ 
can be extended to an Engel structure on $\mti$ 
with a trivial characteristic foliation 
so that it has the given non-negative integer as the minimal twisting number. 
 This implies Theorem~\ref{thmi}-(2). 
 It is also proved that such Engel structures are determined 
by the given contact structures, pairs of Legendrian foliations, 
and non-negative integers. 
 Next, we construct a homotopy between the given Engel structures 
which preserves the trivial characteristic foliation. 
 Then, we apply Theorem~\ref{thm:gol} to show Theorem~\ref{thmi}-(1). 

\subsection{Extension of Engel structures from $M\times \rd I$ to $\mti$.}
\label{sec:ext}
  We show the following Proposition. 
 This implies the second part of Theorem~\ref{thmi}. 
%
%
\begin{prop} \label{prop:determined}
  Let $M$ be a compact orientable $3$--manifold. 
 We suppose that a parallelizable contact structure $\xi$ on $M$, 
a pair $(\F_0,\ \F_1)$ of Legendrian foliations on $(M,\xi)$, 
and a non-negative integer $n\in\nnZ$ are given. 
 Then an Engel structure $\D(\xi,\F_0,\F_1,n)=\D$ 
with a trivial characteristic foliation, 
the induced contact structure $\xi(\D)=\xi$, 
the induced Legendrian foliations $\F(M\times\{i\},\D)=\F_i$, $i=0,1$, 
and the minimal twisting number $\mtw\D=n$,  
is determined up to isotopy by the given $\xi$, $\F_0$, $\F_1$, and $n$. 
\end{prop} 
\begin{proof}
  Since the given contact structure $\xi$ is parallelizable, 
there exists a global framing of $\xi$. 
 Let $(V_0, V_1)$ be a pair of non-vanishing vector fields on $M$ 
which spans $\xi$. 
 We can take this pair 
so that $V_0$ defines the given Legendrian foliation $\F_0$. 
 Then, the Legendrian directions are described by this framing $(V_0,V_1)$. 
 Therefore, we obtain a function $g\co M\to\R$ 
such that the vector field $(V_0\dt\cos g+V_1\dt\sin g)$ generates 
the given Legendrian foliation $\F_1$, 
and that it has its minimum value between $0$ and $\pi$: 
$0<\min_{p\in M}g(p)\le\pi$.
 We define the distribution $\indeng01$ of rank $2$ on $\mti$ 
as the plane field spanned by the vector fields $\uv{t}$ and $V^n$, 
where $t\in I$ is a parameter and 
\begin{equation*}
  V^n:=V_0\dt\cos(t(g+n\pi))+V_1\dt\sin(t(g+n\pi)),\ n\in\nnZ. 
\end{equation*}
 We note that this distribution $\indeng01$ is Engel. 
 It is proved by an easy calculation as follows. 
 The Lie bracket 
\begin{equation*}
  \left[\frac{\rd}{\rd t}, V^n\right]
  =(g+n\pi)\dt\{-V_0\dt\sin(t(g+n\pi))+V_1\dt\cos(t(g+n\pi))\}=:U^n
\end{equation*}
is independent of the distribution $\indeng01$. 
 In addition, the Lie bracket 
\begin{equation*}
  [V^n,U^n]=(g+n\pi)\dt[V_0,V_1]
\end{equation*}
is independent of $\indeng01$ and $U^n$ 
because $(V_0,V_1)$ is a framing of a contact structure, that is, 
a completely non-integrable distribution. 
 Then $\indeng01$ is Engel. 
 The characteristic foliation of this structure is trivial: 
$L\left(\indeng01^2\right)=\spn{\uv t}$. 
 It is nice globally on $\mti$.  
 From the construction, it is clear that the induced contact structure 
by the deprolongation procedure is the given $\xi$. 
 In addition, the induced Legendrian foliations 
on $M\times\{i\},\ (i=0,1)$, is $\F\left(M\times\{i\},\indeng01\right)=\F_i$, 
and the minimal twisting number is $\mtw{\indeng01}=n$. 
 We observe that this structure $\indeng01$ depends 
only on the given $\xi,\F_0,\F_1$, and $n$, up to an isotopy. 
 We show that it is independent of the choices of the framing $(V_0,V_1)$ 
of the given contact structure.  
 Let $(\tld V_0,\tld V_1)$ be another choice of the framing of $\xi$, 
such that $\tld V_0$ defines the Legendrian foliation $\F_0$, 
and that it has the same orientation as the frame $(V_0,V_1)$. 
 Then, there exists a matrix valued function, with a positive determinant, 
\begin{equation*}
  A(p)=
  \begin{pmatrix}
    f_1(p)&0\\ f_2(p)&f_3(p)
  \end{pmatrix}, 
\end{equation*}
which represents a transformation 
between $(\tld V_0,\tld V_1)$ and $(V_0,V_1)$:
\begin{equation*}
  \begin{pmatrix}
    V_0 \\ V_1
  \end{pmatrix}
  =
  \begin{pmatrix}
    f_1(p) & 0 \\ f_2(p) & f_3(p)
  \end{pmatrix}
  \begin{pmatrix}
    \tld V_0 \\ \tld V_1
  \end{pmatrix}. 
\end{equation*}
 We note that $f_1$, $f_3$ are non-vanishing functions with the same sign. 
 Then we have a family of matrix valued functions $A_s(p)$, $s\in [0,1]$, 
between $A_0(p)=A(p)$ and $A_1=\pm\begin{pmatrix}1&0\\0&1\end{pmatrix}$, 
with positive determinants: $|A_s(p)|>0$. 
 We set 
\begin{equation*}
  \begin{pmatrix}
    V_0^s \\  V_1^s
  \end{pmatrix}
  := A_s 
  \begin{pmatrix}
    \tld V_0 \\ \tld V_1
  \end{pmatrix}. 
\end{equation*}
 It is a family of bases of the given contact structure $\xi$, 
between $(V_0^0,V_1^0)=(V_0,V_1)$ and $(V_0^1,V_1^1)=\pm(\tld V_0,\tld V_1)$. 
 We construct a family $\D^s$, $s\in[0,1]$, of Engel structures on $\mti$ 
as above with respect to the bases $(V_0^s,V_1^s)$. 
 We note that we obtain the same Engel structure $\D^1$ 
for both $\pm(\tld V_0,\tld V_1)$. 
 These contact structures coincide at the ends $M\times\rd I$. 
 In fact, the defining vector fields 
$(V^n)_s=V_0^s\dt\cos(t(g_s+n\pi))+V_1^s\dt\sin(t(g_s+n\pi))$ 
are taken so that they determine $\F_i$ on $M\times\{i\}$, $i=0,1$, 
for all $s\in[0,1]$. 
 From the construction, their characteristic foliations are constant: 
$L((\D^s)^2)=\spn{\uv t}$. 
 According to the relative version of Theorem~\ref{thm:gol}, 
this family $\D^s$ is described by a global isotopy. 
 This implies that the induced Engel structure $\indeng01$ is independent 
of the choice of the framing of the given contact structure $\xi$. 
\end{proof}

  By using the development mappings, 
the Engel manifold $(\mti,\indeng01)$ above is regarded 
as a certain subset of the prolongation of $\xi$. 
 We consider the development mapping 
$\dev{\indeng01}\co(\mti,\indeng01)\to(\ptx,\dtx)$. 
 It is an Engel embedding. 
 In this case, the development mapping $\dev{\indeng01}$ is described 
in terms of $V^n$ as follows: 
\begin{equation*}
  \dev{\indeng01}\co(p,t)\mapsto(p,[(V^n)_{(p,t)}]), 
\end{equation*}
where $[(V^n)_{(p,t)}]$ is a point of $\tilde{\pl}(\xi_p)$ 
represented by the line in $\xi_p$ defined by $(V^n)_{(p,t)}$. 
 We set $Q(\xi,\F_0,\F_1,n):=\mathrm{Im}\dev{\indeng01}$ 
and 
\begin{equation*}
  M_n(\F_i):=\Phi_{\indeng01}(M\times\{i\}), \quad i=0,1, 
\end{equation*}
that is, 
\begin{align*}
  &M_n(\F_0)=\{\left(p,[(V_0)_p]\right)\mid p\in M\}, \\
  &M_n(\F_1)=\{\left(p,[V_0\dt\cos(g+n\pi)+V_1\dt\sin(g+n\pi)]\right)
               \mid p\in M\}. 
\end{align*}
 As the given contact structure $\xi$ on $M$ is parallelizable, 
$\ptx$ is identified with $\mtr$ by the relation 
\begin{equation*}
 \Psi_{(V_0,V_1)}\co \ptx\to\mtr, \quad
 (p,[(V_0\dt\cos s+V_1\dt\sin s)_p])\mapsto (p,s), 
\end{equation*}
with respect to the framing $(V_0,V_1)$. 
 By this identification, $M_n(\F_0)$ corresponds 
to the zero-section $M\times\{0\}\sst\mtr$, 
and $M_n(\F_1)$ corresponds to the graph $\{(p,g(p)+n\pi)\in\mtr\mid p\in M\}$ 
of the function $g+n\pi$ on $M$. 

  At the end of this section, we observe a version of the above proposition 
for families of Legendrian foliations and minimal twisting numbers. 
 We suppose that a parallelizable contact structure $\xi$ on $M$, 
a family $(\F_0^s,\F_1^s)$, $s\in[0,1]$, of Legendrian foliations on $(M,\xi)$,
and non-negative integers $n(s)\in\nnZ$ depending on $s\in[0,1]$ are given. 
 We suppose that $n(s)$ changes one by one. 
 Let $(V_0^s,V_1^s)$ be a family of framings of $\xi$, 
such that $V_0^s$ generate $\F_0^s$, $s\in[0,1]$. 
 Then there exists a family $g_s\co M\to\R$ of functions, 
such that $(V_0^s\dt\cos g_s+V_1^s\dt\sin g_s)$ generate $\F_1^s$ 
and $n(s)\pi\le\min_{p\in M}g_s<(n(s)+1)\pi$. 
 In a similar way to the proof of Proposition~\ref{prop:determined}, 
we construct a family $\D(\xi,\F_0^s,\F_1^s,n(s))$ 
of Engel structures on $\mti$ 
as plane fields spanned by $\uv t$ and $(V_0^s\dt\cos g_s+V_1^s\dt\sin g_s)$. 
 By similar arguments to the proof of Proposition~\ref{prop:determined}, 
$\D^s:=\D(\xi,\F_0^s,\F_1^s,n(s))$ is a family of Engel structures on $\mti$ 
with trivial characteristic foliations, 
the induced contact structure $\xi(\D^s)=\xi$, 
the induced Legendrian foliations $\F(M\times\{i\},\D^s)=\F_i^s$, $(i=0,1)$, 
and the minimal twisting numbers $\mtw{\D^s}=n(s)$. 
 The obtained family $\D^s$ of Engel structures is determined up to isotopy 
by the given $\xi$, $\F_0^s$, $\F_1^s$, and $n(s)$, as above. 

\subsection{Proof of Theorem~\ref{thmi}.} \label{sec:pr2}
  In this Section, we prove Theorem~\ref{thmi}-(1). 
 We take $J:=[-1,1]$ as a closed interval instead of $I=[0,1]$, 
for convenience. 
 Let $\D$ be an Engel structure on $\mtj$ 
with the trivial characteristic foliation. 
 This Engel structure has the induced contact structure $\xi=\xi(\D)$, 
the induced Legendrian foliations $\F_i:=\F(M\times\{i\},\D)$, $i=\pm 1$, 
and the minimal twisting number $n=\mtw\D$. 
 We show that the given Engel structure $\D$ is isotopic 
to the Engel structure $\D(\xi,\F_{-1},\F_1,n)$ 
from Proposition~\ref{prop:determined}. 
 This implies that any Engel structure on $\mtj$ 
which has the same induced contact structure $\xi$ on $M$, 
pair of Legendrian foliations $(\F_{-1}$, $\F_1)$ on $M\times\rd J$, 
and minimal twisting number $n$, is isotopic to $\D(\xi,\F_{-1},\F_1,n)$, 
that is, Theorem~\ref{thmi}-(1). 
 We set $W_s:=M\times[-s,s]\sst\mtj$ for $s\in[0,1]$. 
 We note that the induced contact structure on $M$ for $\D\vert_{W_s}$ 
are $\xi=\xi(\D)$ constantly for all $s\in[0,1]$. 
 We set $\F_t:=\F(M\times\{t\},\D),\ t\in[-1,1]$. 
 Identifying $M\times\{t\}\in\mtj$ with $M$ by the standard projection, 
$\F_t$ can be regarded as a family of Legendrian foliations on $(M,\xi)$. 
 Let $n(s):=\mtw{\D\vert_{W_s}}$ be the minimal twisting number of $\D$ 
restricted to $W_s$. 
 We note that $n(1)=\mtw\D$, $n(\eps)=0$ for sufficiently small $\eps>0$, 
and that it is monotone with respect to $s\in(0,1]$. 
 According to Proposition~\ref{prop:determined} and its version for families, 
we obtain a family $\tld\D_s,\ s\in[0,1]$, of Engel structures on $\mti$ 
which have 
$$\xi(\tld\D_s)=\xi,\quad \F(M\times\{0\},\tld\D_s)=\F_{-s},\quad 
  \F(M\times\{1\},\tld\D_s)=\F_s,$$ 
$$\text{and}\quad \mtw{\tld\D_s}=n(s).$$ 
Applying a family of diffeomorphisms $\psi_s\co\mti\to W_s=M\times[-s,s]$, 
preserving fibres, 
we obtain a family $\D_s$ of Engel structures on each $W_s$ 
with the trivial characteristic foliation, 
which is endowed with the given informations: 
\begin{equation*}
  \xi(\D_s)=\xi,\ \F(M\times\{\pm s\},\D_s)=\F_{\pm s},\ 
  \text{and}\ \mtw{\D_s}=n(s). 
\end{equation*}
 Now, we consider the following family $\E_s,\ s\in[\eps,1]$, 
of Engel structures on $\mtj$, where $\eps>0$ is sufficiently small, 
\begin{equation*}
  \E_s:=
  \begin{cases}
    \D \quad   &\text{on}\ (\mtj)\setminus W_s, \\
    \D_s \quad &\text{on}\ W_s. 
  \end{cases}
\end{equation*}
 We note that Theorem~\ref{thm:mont} guarantees the pasting of Engel structures
along $M\times\{\pm s\}$. 
 This family is a path, relative to the ends $M\times\rd J$, 
between $\E_{\eps}$ and $\E_1=\D_1$. 
 The structure $\D_1=\D(\xi,\F_0,\F_1,n)$ is determined, up to isotopy, 
only by the given contact structure, pair of Legendrian foliations, 
and non-negative integer. 
 According to Theorem~\ref{thm:mont}, the germ of Engel structure 
along $M\times\{0\}$ is determined only by the induced Legendrian foliation 
$\F(M\times\{0\},\D)$ on $M\times\{0\}$. 
 Therefore, the structure $\E_{\eps}$ is isotopic to the given structure $\D$. 
 We note that the family $\E_s$ above preserves 
the trivial characteristic foliation: $\cL({\E_s}^2)\equiv\spn{\uv t}$, 
so we can apply the relative version of Theorem~\ref{thm:gol}. 
 Consequently, the given Engel structure is isotopic, 
relative to the ends $M\times\rd I$, 
to the Engel structure 
\begin{equation*}
  \D(\xi(\D),\F(M\times\{-1\}),\F(M\times\{1\}),\mtw\D), 
\end{equation*}
which depends only on the induced contact structure 
and Legendrian foliations on both ends, and the minimal twisting number. 
\endproof

  Let us conclude this paper by an open problem. 
 The results of this paper need some conditions. 
 It is natural and might be interesting to consider more general cases. 
%
%
\begin{qtn}
  Let $\xi_0$, $\xi_1$ be contact structures on a $3$--manifold $M$, 
which are not isomorphic to each other. 
 Is there any Engel structure $\D$ on a $4$--manifold $\mti$, 
whose characteristic foliation $\cL(\D^2)$ is transverse 
to both $M\times\{0\}$, $M\times\{1\}$, 
and which induce the given contact structures $\xi_0$, $\xi_1$ 
on $M\times\{0\}$, $M\times\{1\}$, respectively? 
\end{qtn} 

  We note that there must be leaves of characteristic foliation 
which do not intersect both ends $M\times\{0\}$, $M\times\{1\}$, 
if such an Engel structure exists. 

%
%

\Addresses\recd


\begin{thebibliography}
\bibitem[A]{adachi} J.~Adachi, 
  \textit{Germs of Engel structures along $3$--manifolds}, preprint. 

\bibitem[BCG3]{bcg} R.~L.~Bryant, S.~S.~Chern, R.~B.~Gardner, 
  H.~L.~Goldschmidt, P.~A.~Griffiths, 
  \textit{Exterior differential systems}, 
  Mathematical Sciences Research Institute Publications, \textbf{18}.
  Springer-Verlag, New York, 1991.

\bibitem[C]{cartan} E.~Cartan, 
  \textit{Sur quelques quadratures dont l'element differentiel 
          contient des fonctions arbitraires}, 
  Bull.\ Soc.\ Math.\ France \textbf{29} (1901), 118--130. 

\bibitem[El1]{classot} Ya.~Eliashberg, 
  \textit{Classification of overtwisted contact structures on $3$--manifolds}, 
  Invent.\ Math.\ \textbf{98} (1989), 623--637. 

\bibitem[El2]{20years} Ya.~Eliashberg, 
  \textit{Contact $3$--manifolds twenty years since J. Martinet's work}, 
  Ann.\ Inst.\ Fourier (Grenoble) \textbf{42} (1992), no. 1-2, 165--192. 

\bibitem[Et]{etnyre} J.~B.~Etnyre, 
  \textit{Tight contact structures on lens space}, 
  Commun.\ Contemp.\ Math.\ \textbf{2} (2000), no.\ 4, 559--577.

\bibitem[Ge]{ger} V.~Gershkovich, 
  \textit{Exotic Engel structures on $\R^4$}, 
  Russian J. Math.\ Phys.\ \textbf{3} (1995), no.\ 2, 207--226.

\bibitem[Go]{gol} A.~Golubev, 
  \textit{On the global stability of maximally nonholonomic two-plane fields 
          in four dimensions}, 
  Internat.\ Math.\ Res.\ Notices 1997, no.\textbf{11}, 523--529. 

\bibitem[Gi]{giroux} E.~Giroux,
  \textit{Structures de contact en dimension trois 
  et bifurcations des feuilletages de surfaces}, 
  Invent.\ Math.\ \textbf{141} (2000), no.\ 3, 615--689.

\bibitem[H]{honda} K.~Honda, 
  \textit{On the classification of tight contact structures I}, 
  Geom.\ Topol.\ \textbf{4} (2000), 309--368. 

\bibitem[K]{kanda} Y.~Kanda, 
  \textit{The classifications of contact structures on $3$--torus}, 
  Comm.\ Anal.\ Geom.\ \textbf{5} (1997), 413--438.

\bibitem[M1]{mont1} R.~Montgomery, 
  \textit{Generic distributions and Lie algebras of vector fields}, 
  J. Differential Equations \textbf{103} (1993), no. 2, 387--393. 

\bibitem[M2]{mont2} R.~Montgomery, 
  \textit{Engel deformations and contact structures}, 
  Northern California Symplectic Geometry Seminar, 
  Amer.\ Math.\ Soc.\ Transl.\ Ser.\ 2, \textbf{196}, 103--117, 
  Amer.\ Math.\ Soc., Providence, RI, 1999. 

\bibitem[VG]{vege} A.~Vershik and V.~Gershkovich, 
  \textit{Nonholonomic dynamical systems. Geometry of distributions 
          and variational problems}, 
  Encyclopaedia Math.\ Sci.\ \textbf{16} (1994), 1--81. 
\end{thebibliography}
\end{document}
